\documentclass[12pt]{amsart}
\usepackage{hyperref}
\usepackage{amssymb}
\usepackage{amsmath}
\usepackage{amsthm}
\usepackage{graphicx}
\usepackage{color}
\usepackage{listings}
\usepackage{subfigure,multirow,booktabs,url,breakurl,enumerate}
\usepackage[table]{xcolor}
\newcommand\shade{0.9}
\definecolor{lightgrey}{rgb}{\shade,\shade,\shade}

\newtheorem{theorem}{Theorem}
\newtheorem{corollary}{Corollary}
\newtheorem{lemma}{Lemma}

\newtheorem{conjecture}{Conjecture}

\newcommand{\lref}[1]{Lemma~\ref{lem:#1}}
\newcommand{\cref}[1]{Corollary~\ref{cor:#1}}
\newcommand{\conref}[1]{Conjecture~\ref{conj:#1}}
\newcommand{\eref}[1]{Equation~(\ref{eq:#1})}
\newcommand{\sref}[1]{Section~\ref{sec:#1}}
\newcommand{\tref}[1]{Theorem~\ref{thm:#1}}

\newcommand{\tabref}[1]{Table~\ref{tab:#1}}
\newcommand{\fref}[1]{Fig.~\ref{fig:#1}}
\newcommand{\lstref}[1]{Listing~\ref{lst:#1}}
\newcommand\set[1]{\{ #1\}}
\newcommand\flr[2]{\left\lfloor \frac{#1}{#2} \right\rfloor}

\newcommand\ldown[1]{l_{#1}^\downarrow}
\newcommand\lup[1]{l_{#1}^\uparrow}
\newcommand\rdown[1]{r_{#1}^\downarrow}
\newcommand\rup[1]{r_{#1}^\uparrow}

\newcommand\tright[1]{t_{#1}^\rightarrow}
\newcommand\tleft[1]{t_{#1}^\leftarrow}
\newcommand\bright[1]{b_{#1}^\rightarrow}
\newcommand\bleft[1]{b_{#1}^\leftarrow}

\newcommand\ssum[1]{\sum_{ \substack{#1}}}
\newcommand{\Tn}[1]{\mathbf T_{#1}}
\newcommand{\vd}{V}
\newcommand{\hd}{H}
\newcommand{\vh}{V\hspace{-5pt}H}
\newcommand{\round}[2]{\left\lfloor
    \frac{#1}{#2} + \frac{1}{2} \right\rfloor}
\newcommand{\op}{Od}
\newcommand{\I}[1]{{[\![#1]\!]}}            % Iversonian

\newcommand{\fd}{`$\diagup$'}
\newcommand{\bd}{`$\diagdown$'}

\begin{document}
\pagestyle{headings}

\lstset{language=Python,%
  basicstyle=\small\ttfamily,%
  keywordstyle=\color{blue}\bfseries,%
  commentstyle=\color{gray},%
  stringstyle=\color{purple}\ttfamily,%
  showstringspaces=false,%
%  numbers=left,%
%  numberstyle=\small,%
%  stepnumber=1%
}

%NON ELS class
\title[$n\times n$ tatami coverings with $v$ vertical
dominoes]{Enumerating maximal tatami mat coverings of square grids
  with $v$ vertical dominoes}%
\author[Erickson]{Alejandro Erickson}
% \address{Department of Computer Science, University of Victoria, V8W
%   3P6, Canada}

  \author[Ruskey]{Frank Ruskey}
  \address{Department of Computer Science, University of Victoria, V8W
    3P6, Canada}

  \maketitle

  %%%%%%%%%%%%%%%%%%%%% ELS class only
% \begin{frontmatter}
%   \title{Enumerating tatami mat arrangements of square grids with $v$ vertical dominoes}
%   \author{Alejandro Erickson\corref{cor}}
%   \cortext[cor1]{Correpsonding author}
%   \address{Department of Computer Science, University of Victoria, V8W
%     3P6, Canada}
%   \ead{ate@uvic.ca}
%   \ead[url]{http://alejandroerickson.com}

%   \author{Frank Ruskey}
%   \address{Department of Computer Science, University of Victoria, V8W
%     3P6, Canada}
%   \ead{ruskey@cs.uvic.ca}
%   \ead[url]{http://webhome.cs.uvic.ca/~ruskey/}

\begin{abstract}
  We enumerate a certain class of monomino-domino coverings of square
  grids, which conform to the \emph{tatami} restriction; no four tiles meet.

  Let $\Tn{n}$ be the set of monomino-domino tatami coverings of the
  $n\times n$ grid with the maximum number, $n$, of monominoes,
  oriented so that they have a monomino in each of the top left and
  top right corners.  We give an algorithm for exhaustively generating
  the coverings in $\Tn{n}$ with exactly $v$ vertical dominoes in
  constant amortized time, and an explicit formula for counting them.
  The polynomial that generates these counts has the factorisation
  \begin{align*}
    P_n(z)\prod_{j\ge 1} S_{\left\lfloor \frac{n-2}{2^j}
      \right\rfloor}(z),
  \end{align*}
  where $S_n(z) = \prod_{i=1}^{n} (1 + z^i)$, and $P_n(z)$ is an
  irreducible polynomial, at least for ${1 < n < 200}$.  We present some
  compelling properties and conjectures about $P_n(z)$.  For example
  $P_n(1) = n2^{\nu(n-2)-1}$ for all $n \ge 2$, where $\nu(n)$ is the
  number of $1$s in the binary representation of $n$ and $\deg(P_n(z))
  = \sum_{k=1}^{n-2} \op(k)$, where $\op(k)$ is the largest odd
  divisor of $k$.
 \end{abstract}

 %%%%%%%%%%%%%%%%%%%%%%%%%%%  ELS class only
%  \begin{keyword}
%    %mostly copied over from TNNM for JDA.
%    tatami\sep monomino-domino covering\sep monomino\sep domino\sep
%    covering\sep polyomino
%    \MSC[2010] 05B45 \sep \MSC[2010] 05B50
%    \end{keyword}
% \end{frontmatter}

\section{Introduction}
The counting of domino coverings, together with its extension to
counting perfect matchings in (planar) graphs, is a classic area of
enumerative combinatorics and theoretical computer science.  Less
attention has been paid, however, to problems where the local
interactions of the dominoes are restricted in some fashion.  Perhaps
the most natural such restriction is the ``tatami'' condition, defined
below.  The tatami condition is quite restrictive: for example, the
$10\times 13$ grid cannot be covered with dominoes and also satisfy
the tatami condition.  In this paper we restrict our attention to
square grids, and explore in some detail the enumeration of certain
extremal configurations.

Tatami mats are a traditional Japanese floor covering whose dimensions
are approximately $1\textup{m} \times 1\textup{m}$ or
$1\textup{m}\times
2\textup{m}$. %following 3.27 of chicago man style 16th ed.
In certain arrangements, no four tatami mats may meet.  Such an
arrangement has a preferable structure which is discussed in
\cite{EricksonRuskeySchurch2011} and \cite{EricksonSchurch2012}.

A \emph{tatami covering} is an arrangement of $1\times 1$ monominoes,
$1\times 2$ horizontal dominoes, and $2\times 1$ vertical dominoes, in
which no four tiles meet.  The present discussion is about tatami
coverings of the $n\times n$ grid with exactly $n$ monominoes and $v$
vertical dominoes.  On the basis of some computer investigations, Don
Knuth discovered that the generating polynomial for small tatami
coverings of this type, with respect to the number of vertical
dominoes they contain, is a product of cyclotomic polynomials and a
mainly mysterious, irreducible polynomial (private communication,
December 2010). Knuth's discovery and our own observations motivated
Conjecture 4 in \cite{EricksonRuskeySchurch2011}, which is presented
here as \eref{factorisation}.  In this paper we generalize and prove
Knuth's cyclotomic factors, and determine some important properties of
the mysterious polynomial.

We prove in \cite{EricksonRuskeySchurch2011} that all of the $n\times
n$ coverings with $n$ monominoes can be rotated so that monominoes
appear in each of the top two corners of the grid, so we let $\Tn{n}$
be the set of these.  Let $H(n,k)$ be the number of coverings in
$\Tn{n}$ with exactly $k$ horizontal dominoes, and let $V(n,k)$ be the
number with exactly $k$ vertical dominoes.  Let $S_n(z) =
\prod_{i=1}^n (1+z^i)$.  We prove that the polynomial
\begin{align}
  %
  % Verified against HnkVnk.sage on June 13, 2012
  %
  \label{eq:tnzintro}
  \vh_n(z):=2\sum_{i=1}^{\flr{n-1}{2}}
  S_{n-i-2}(z)S_{i-1}(z)z^{n-i-1} +
  \left(S_{\flr{n-2}{2}}(z)\right)^2,
\end{align}
is equal to $\sum_{k\ge 0} H(n,k)z^k$ for odd $n$ and is equal to
$\sum_{k\ge 0} V(n,k)z^k$ for even $n$.  Knuth's observation
generalizes to
\begin{align}
  \label{eq:factorisation}
  \vh_n(z) = P_n(z) \prod_{j\ge 1} S_{\left\lfloor \frac{n-2}{2^j}
    \right\rfloor}(z),
\end{align}
where $P_n(z)$ is the ``mysterious'' polynomial.  We prove
\eref{factorisation} in \tref{cyclotomic}.

The remaining factors of \eref{factorisation} are of the form
$S_k(z)$, where $k$ is a binary right shift of $n-2$, and the complete
factorisation of these is known in general.  The \emph{$i$th
  cyclotomic polynomial}, $\Phi_i(z)$, is defined as
$\prod_{\omega\in\Omega} (z-\omega)$, where $\Omega$ is the set of
$i$th primitive roots of unity.  Lemma 5 in
\cite{EricksonRuskeySchurch2011} shows that $S_k(z)$ is a certain
product of cyclotomic polynomials, which are known to be irreducible,
and thus $\vh_n(z)$ can apparently be factored completely as
\begin{align}
  % Proofed June 13 2012
  \label{eq:completeF}
  \vh_n(z) = P_n(z) \prod_{j\ge 1}  \Phi_{2j} (z)^{\flr{n-2}{2j}}.
\end{align}

We have verified the irreducibility of $P_n(z)$ for ${1< n < 200}$ (the
degree of $P_{199}(z)$ is $13022$ and its largest coefficient has $55$
digits), and thus we hope that \eref{completeF} is the complete
factorisation of $\vh_n(z)$ for all $n\ge 2$.

The class $P_n(z)$ of polynomials has some compelling properties, some
of which are theory, others empirical.  For example, we observe in
\conref{Pcon} that the alternating sums of $P_n(z)$ are the
coefficients of the ordinary generating function
\begin{align*}
  \sum_{n\ge 2} P_n(-1)z^{n-2} =
  \frac{(1+z)(1-2z)}{(1-2z^2)\sqrt{1-4z^2}},
\end{align*}
for ${1 < n < 200}$.  If the conjecture is true, then $P_{2(n+1)}(-1) =
\binom{2n}{n}$. Furthermore $P_n(-1)$ is equal to the sum of the
absolute values of the coefficients of $P_n(z)$, only for $n\ge 20$.
This second fact is surprising, considering the way $P_n(z)$ is
derived -- why $n\ge 20$?

The complex roots of $P_n(z)$ appear to cluster neatly around the unit
circle, and form convergent sequences as $n\longrightarrow \infty$.
They are plotted in \fref{complexPlotP} for odd $n$; for even $n$, the
plot has a similar look.

Theoretical progress on $P_n(z)$ comprises \tref{degpn} and
\tref{nu}.  The former states that $\deg(P_n(z)) = \sum_{k=1}^{n-2}
\op(k)$, where $\op(n)$ is the largest odd divisor of $n$. We prove in
\tref{nu} that for all $n\ge 2$, the sum of the coefficients of
$P_n(z)$ is equal to $n2^{\nu(n-2)-1}$, where $\nu(n)$ is the number
of $1$-bits in the binary representation of $n$.

Our technique for finding $\vh_n(z)$ employs an operation which
preserves the tatami condition, called the \emph{diagonal flip},
defined in \cite{EricksonSchurch2012}.  The added observation that a
diagonal flip changes the orientation of some dominoes, enables us to
further exploit it.  The crux of the argument uses the partition of
$\Tn{n}$, from Theorem 2 of \cite{EricksonSchurch2012}, which reveals
diagonal flips each with $1,2,\ldots, k$ dominoes, respectively, that
can be flipped independently.  We use this to express $\vh_n(z)$ in
terms of $S_k(z)$, the generating polynomial for the number of subsets
of $\{1,2,\ldots,k\}$ whose elements sum to $i$.

The formula for the $v$th coefficient of $\vh_n(z)$ translates into an
algorithm for generating all possible $n\times n$ tatami coverings
with $v$ vertical (or horizontal) dominoes, given that we have one for
generating all $k$-sum subsets of the $n$ set, for $k,n\ge 0$.  We
employ an algorithm from \cite{BaronaigienRuskey1993} to generate our
coverings in constant amortized time.

\subsection{Overview}
The structure of square tatami coverings with the maximum number of
monominoes is summarized in \sref{structure} --- see reference
\cite{EricksonRuskeySchurch2011} for a complete proof.  In
\sref{representation} we describe a representation for the coverings
in $\Tn{n}$ using strings over a ternary alphabet.  Each symbol
represents a monomino which can either be flipped in exactly one of
two diagonals, or unflipped.  This representation is essentially the
same as the one in \cite{EricksonSchurch2012}, but our revision of the
indices greatly simplifies the counting of vertical and horizontal
dominoes.

Our main result is \tref{vdominoes}, the generating polynomial
$\vh_n(z)$, but the technical part is proving the formula for its
coefficients, which is presented in \lref{vdominoes}.  We discuss the
coefficients of $\vh_n(z)$ and derive the self-reciprocal generating
polynomial required to partition the coverings in $\Tn{n}$ and their
four distinct rotations according to their numbers of vertical
dominos.

In \sref{combalg} we apply the proof of \lref{vdominoes} to generate
$\Tn{n}$ in constant amortized time by adapting an algorithm given in
\cite{BaronaigienRuskey1993}.

In \sref{pn} we depart from the geometrical interpretation of
$\vh_n(z)$ and prove \tref{cyclotomic}, the factorisation in
\eref{factorisation}, and the remainder of the section focuses on
properties of the factor $P_n(z)$, including \tref{degpn} and
\tref{nu}, introduced above.

\section{Structure of a tatami covering}
% Proofed June 13, 2012
\label{sec:structure}
An \emph{edge} refers to the edge of a tile, while a \emph{boundary}
refers to the outer boundary of a covering.

The structure of coverings in $\Tn{n}$ is characterized in Lemma 3 and
Corollary 2 of \cite{EricksonRuskeySchurch2011}.  Corollary 2 states
that an $n\times n$ covering with $n$ monominoes has monominoes in
exactly two of its corners, which must share a boundary.  This ensures
that each covering of $\Tn{n}$ has distinct rotations through
$0,\pi/2,\pi,$ and $3\pi/2$ radians, and that rotations of distinct
coverings in $\Tn{n}$ are distinct from each other.

If $T\in \Tn{n}$, then a \emph{diagonal}, $D$, of $T$ is a contiguous
sequence of like-aligned dominoes whose centers lie on a line with
slope $1$ or $-1$.  The sequence must begin with a domino with its
long edge on the boundary; the final domino will share an edge with a
monomino which is also considered to be part of the diagonal.  In
\fref{diag} there are 12 diagonals, only one of which contains
horizontal dominoes.

A \emph{diagonal flip} of $D$ consists of removing it from $T$,
reflecting horizontally, rotating by $\frac{\pi}{2}$ radians, and
placing it back onto the grid squares that were vacated.

There are three things to note about the diagonal flip:
\begin{itemize}
\item a flipped diagonal is a diagonal;
\item the operation preserves the tatami restriction; and,
\item it changes the orientation of the dominoes that it contains, and
  maps the monomino to the other extreme of the diagonal.
\end{itemize}

The \emph{running bond}, or simply \emph{bond}, is a rotation of the
basic brick laying pattern, in which all dominoes have the same
orientation.  The restriction that coverings in $\Tn{n}$ have
monominoes in their upper corners, implies that exactly one bond
pattern is possible for each $n$.  When $n$ is even, the bond consists
of horizontal dominoes, with monominoes along the left and right
boundaries, and when $n$ is odd, the bond consists of vertical
dominoes, with monominoes along the top and bottom boundaries (see
\fref{bond}).

Lemma 3 in \cite{EricksonRuskeySchurch2011} shows that every covering in
$\Tn{n}$ can be produced from the running bond, via a finite sequence
of diagonal flips in which each monomino is moved at most once, and the
top corner monominoes are not moved at all. In addition, any such
sequence of flips results in an element of $\Tn{n}$ (see \fref{flips}).

\begin{figure}[h]
  \centering \subfigure[ ]{%
    \includegraphics{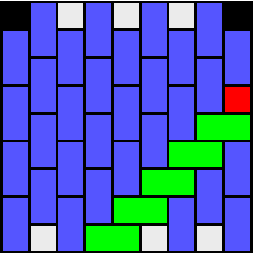}%
    \label{fig:diag}}%
  \hspace{0.5in}%
  \centering \subfigure[ ]{%
    \includegraphics{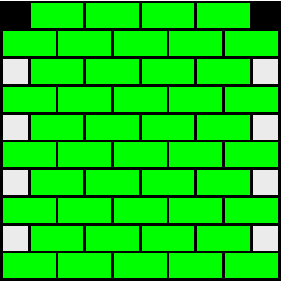}%
    \label{fig:bond}}%
  \hspace{0.5in}%
  \subfigure[ ]{%
    \includegraphics{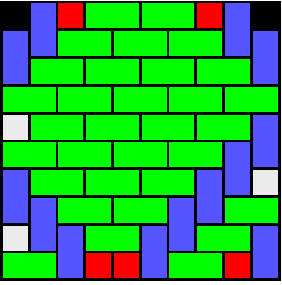}%
    \label{fig:flips}}%
  \caption{Examples of coverings in $\Tn{n}$. \emph{(a)} A diagonal
    flip in a $9\times 9$ covering. \emph{(b)} The horizontal running
    bond for $n=10$. \emph{(c)} A sequence of diagonal flips results
    in a covering in $\Tn{10}$.  Flipped monominoes are coloured red.}
  \label{fig:bondflips}
\end{figure}

From this perspective, we can look upon the original position of a
monomino as its place in the running bond, and then describe it as
\emph{flipped} in a given direction if it has moved from its original
position in $T$; otherwise it is \emph{unflipped}.

% Superficially, we have departed from the tatami condition, that no
% four tiles meet.  Let us describe $\Tn{n}$ a different way, for even
% $n$.  Begin with an $n\times n$ chess board (with a white square in
% the bottom right), and place a bishop on every white square along
% the left boundary, and on every black square along the right
% boundary, except the top right.  That is, each bishop represents a
% monomino that can be flipped.  Now make at most one move of maximal
% length for each bishop, such that no two moves cross paths.  Each
% move corresponds to a diagonal flip, and the restriction that moves
% do not cross imply that this is an element of $\Tn{n}$.  The total
% number of squares travelled by bishops, using the manhattan metric,
% is twice the number of dominoes that change orientation in the
% corresponding covering.

\section{Representing coverings as a ternary string}
\label{sec:representation}
We describe a ternary string representation for $n\times n$ coverings
with $n$ monominoes.  Recall that each monomino, besides the two
corner monominoes, is in exactly two diagonals in the running bond,
and in a given covering a monomino is flipped in one of these
diagonals, or it is unflipped.  A ternary symbol for each monomino
indicates which of the three possible states it assumes.  Each
covering is described by a unique string of these ternary symbols,
represented in the same order as the following indexed labelling (see
caption at \fref{ternary}).

Monominoes and their diagonals are labelled as shown in
\fref{labelling}, such that the index, $i$, of a monomino is equal to
the length of one of its diagonals, and $n-i-1$ is the length of the
other.  This relationship between diagonal length and index is helpful
in \lref{vdominoes}.

The ternary string representing the $10\times 10$ covering in
\fref{ternary} is $s=(0,1,-1,0,0,1,-1)$, where $s_i=1$ if the
$i^\textup{th}$ monomino is flipped upward, $s_i=-1$ if it is flipped
downward, and $s_i=0$ if it is unflipped.

\begin{figure}[h]
  % Proofed June 13, 2012
  \centering \subfigure[]{
    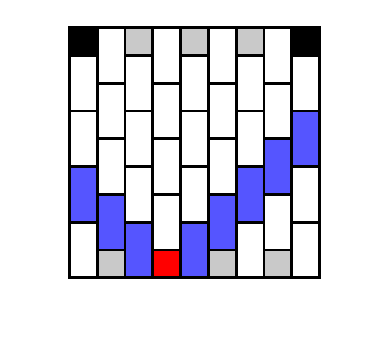%
    \label{fig:oddlabelling}}%
  \subfigure[]{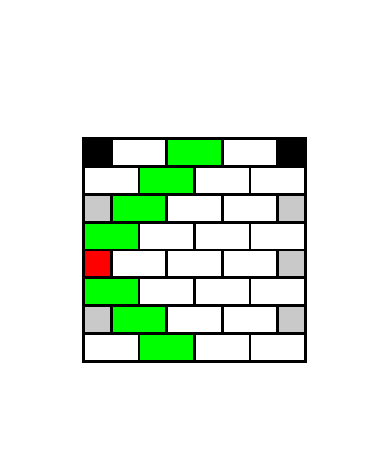%
    \label{fig:evenlabelling}}%
  \subfigure[]{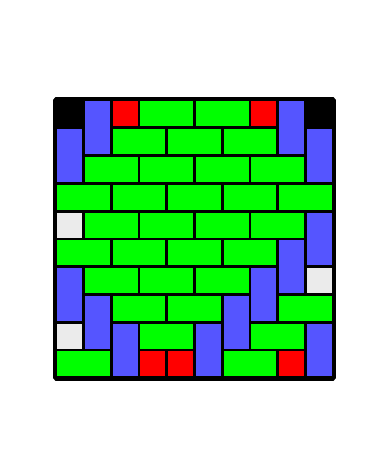%
    \label{fig:ternary}}%
  \caption{Labelling for $\Tn{n}$. \emph{(a)} For odd $n$, monominoes
    are labelled $t_i$ and $b_i$.  The distances from $t_i$ and $b_i$
    to the left boundary are both $i$. \emph{(b)} For even $n$,
    monominoes are labelled $l_i$ and $r_i$.  The distances from $l_i$
    to the bottom boundary, and from $r_i$ to the top boundary, are
    both $i$.  \emph{(c)} The covering, $(0,1,-1,0,0,1,-1)$.}
  \label{fig:labelling}
\end{figure}

We use $\lup{i}, \ldown{i},
\rup{i},\rdown{i},\tright{i},\tleft{i},\bright{i},\bleft{i}$ to denote
the diagonals that the monominoes $l_i,r_i,t_i,b_i$ can be flipped on.
Naturally, $l_i$ and $r_i$ can only be (diagonally) flipped up or
down, whilst $t_i$ and $b_i$ can only be flipped left or right.

Let $d_n(a)$ be the number of dominoes in the diagonal $a$, also called
the length or size of the diagonal.  It is a function of the index and
direction of $a$:

\begin{align*}
  % Proofed June 13, 2012
  d_n(a)=\left\{\begin{aligned} i&,& \text{if }& a\in \set{\ldown{i},
      \rup{i}, \tleft{i}ñ, \bleft{i}};\\
      n-i-1&,& \text{if }& a\in\set{\lup{i},\rdown{i},\tright{i},\bright{i}}.\\
    \end{aligned}\right.
\end{align*}

Flipped diagonals which intersect are called \emph{conflicting}, and can occur as one of two types (see \fref{conflicts}).

\begin{figure}[h]
  \centering \subfigure[]{%
    \includegraphics{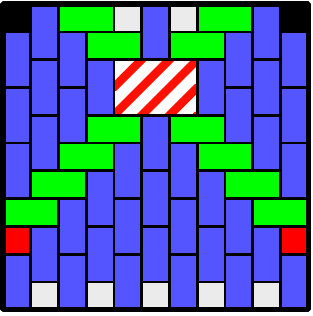}%
    \label{fig:type1}}%
  \hspace{0.5in}%
  \centering \subfigure[]{%
    \includegraphics{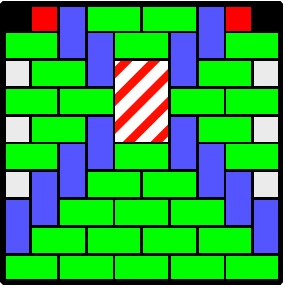}%
    \label{fig:type2}}%
  \hspace{0.5in}%

  \caption{Example of, \emph{(a)}, Type 1 conflict, and, \emph{(b)},
    Type 2 conflict.}
  \label{fig:conflicts}
\end{figure}

\begin{description}
\item[Type 1] A pair of diagonals with monominoes originating on the
  same boundary are flipped toward one another (e.g. $(\tright{i},
  \tleft{j})$ for some $i<j$).
\item[Type 2] A pair of diagonals with monominoes originating on
  opposite boundaries are flipped in the same direction
  (e.g. $(\lup{i},\rup{j})$) and their combined length is at least
  $n$ (see \tabref{conflicts}).
\end{description}

\begin{table}
\begin{tabular}[h]{ll}
    % Proofed June 13, 2012
  \toprule
  Pair & Type 2 $\iff$\\
  \hline
  $\ldown{i},\rdown{j}$ &$ j \le i-1$\\%$n \le i + n-j-1$\\
  $\lup{i}, \rup{j}$& $ i\le j-1$ \\ %n \le n-i-1 + j
  $\tleft{i}, \bleft{j}$& $ n\le j+i$ \\ %n \le i + j
  $\tright{i}, \bright{j}$& $i+j \le n-2  $ \\%$ n-i-1+n-j-1\ge n$
  \bottomrule
\end{tabular}
\vspace{0.05in}
  \caption{Conditions for Type 2 conflicts.}
  \label{tab:conflicts}
\end{table}

Lastly, if $a$ is a diagonal containing a given monomino, let $\bar a$
be the monomino's other diagonal.

\subsection{A partition of $\Tn{n}$}

Let $\Tn{n}(a)\subseteq \Tn{n}$, where $a$ is a diagonal such that
$d_n(a) \ge d_n(\bar a)$, be defined as the collection of coverings in
$\Tn{n}$ in which $a$ is the \emph{longest flipped diagonal}; for each
flipped diagonal $b$, distinct from $a$, we have $d_n(b)<d_n(a)$.

Let $\Tn{n}(\varnothing)$ be the set of coverings in which no monomino is
flipped on its longest diagonal.  Note the distinction between a
monomino flipped on its longest diagonal, and the longest flipped
diagonal in the whole covering.

The sets $\Tn{n}(\varnothing)$ and $\Tn{n}(a)$, for each diagonal $a$
defined above, are a partition of $\Tn{n}$, and the allowable diagonal
flips of each subset can be applied independently of the other flips,
by Theorem 2 in \cite{EricksonSchurch2012}.

\section{Enumeration}
\label{sec:enumeration}
Let $S(s,k)$ be the number of subsets of $\{1,2,\ldots, s\}$ whose sum
is $k$.  The number of coverings with $k$ vertical (or horizontal)
dominoes is expressible in terms of this function by making
independent flips of diagonals whose lengths are some subset
$\{1,2,\ldots, s\}$.  We identify these sets of diagonals in the proof
of \lref{vdominoes}.

\begin{lemma}
  \label{lem:vdominoes}
  Let $\vd(n,k)$ and $\hd(n,k)$ be the number of coverings in $\Tn{n}$
  with exactly $k$ vertical and horizontal dominoes, respectively.  If
  $n$ is even, then $\vd(n,k)$ is equal to
  \begin{subequations}
    \begin{align}
      % Proofed/Verified against HnkVnk.sage June 13, 2012
      \label{eq:tnka}%
      \vh(n,k) :=& 2\sum_{i=1}^{\flr{n-1}{2}} \left( \ssum{k_1+k_2 = \\k
          -(n-i-1)}
        S(n-i-2,k_1)S(i-1,k_2) \right)\\
      \label{eq:tnkb}%
      + & \sum_{k_1+k_2 = k}
      S\left(\flr{n-2}{2},k_1\right) S\left(\flr{n-2}{2},k_2\right).
    \end{align}
  \end{subequations}
  When $n$ is odd, $\vh(n,k)$ is equal to $\hd(n,k)$.
\end{lemma}
\proof Each outer sum term of (\ref{eq:tnka}) adds the coverings for
$\Tn{n}(a)$, for some diagonal $a$, and the term (\ref{eq:tnkb})
counts those in $\Tn{n}(\varnothing)$.

\textbf{Case $n$ even:} The trivial covering in $\Tn{n}$ consists only
of horizontal dominoes, and flipping the diagonal $a$ contributes
$d_n(a)$ vertical dominoes.  Diagonals $\lup{i}$ and $\rup{i}$ have even
length, for all $i$, while $\ldown{i}$ and $\rdown{i}$ have odd
length.  We use this fact to find sets of diagonals which have lengths
$1,2,\ldots, s$, for some $s\in \mathbb N$, by combining allowable
diagonals in opposite corners, for each $\Tn{n}(a)$.
\tabref{tnnnCases} shows the lengths of the longest allowable diagonals in each
corner for each $\Tn{n}(a)$, and from this we can find the required
sets of diagonals.  For example, the allowable diagonals in
$\Tn{n}(\lup{i})$ are shown in \fref{evenExample} (for $(n,i)=(18,5)$)
and their respective lengths are
\begin{center}
%VERIFIED FOR CORRECTNESS AND CONSISTENCY WITH TABLE 2 BY ALEJANDRO ON JAN 12, 2012
  \begin{tabular}{ll}
      % Proofed June 13, 2012
  $    \ldown{1},\ldown{3}, \ldots, \ldown{i-2} $&$  1,3, \ldots, i-2,$\\
  $    \lup{i+2},\lup{i+4}, \ldots, \lup{n-3}$&$  n-i-3, n-i-5, \ldots, 2,$\\
  $    \rdown{i+1}, \rdown{i+3}, \ldots,\rdown{n-2}$&$ n-i-2, n-i-4,
  \ldots,1,$\\
  $    \rup{2}, \rup{4}, \ldots,\rup{i-1}$&$  2, 4, \ldots,i-1.$
\end{tabular}
\end{center}

\begin{figure}[h]
  \centering
  \subfigure[]{%
    {\scriptsize \def\svgwidth{0.48\textwidth}%
      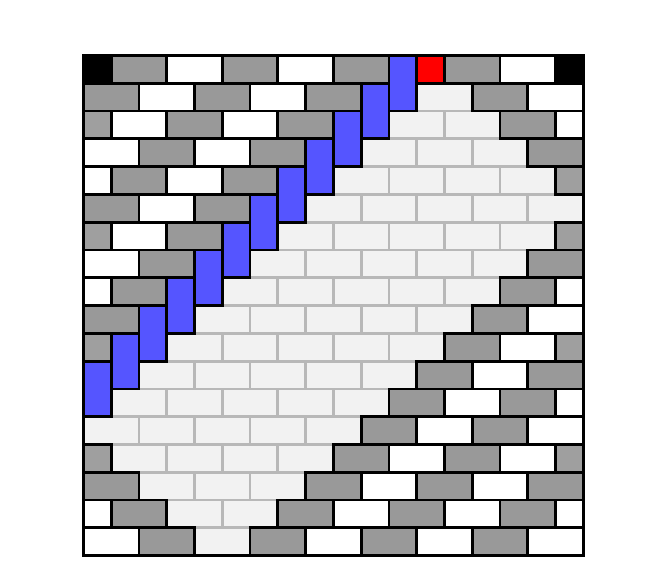%
      }
    \label{fig:evenExample}} %
  \hspace{0.00\textwidth}%
  \subfigure[]{%
    {\scriptsize \def\svgwidth{0.48\textwidth}%
      \def\svgwidth{0.48\textwidth}%
      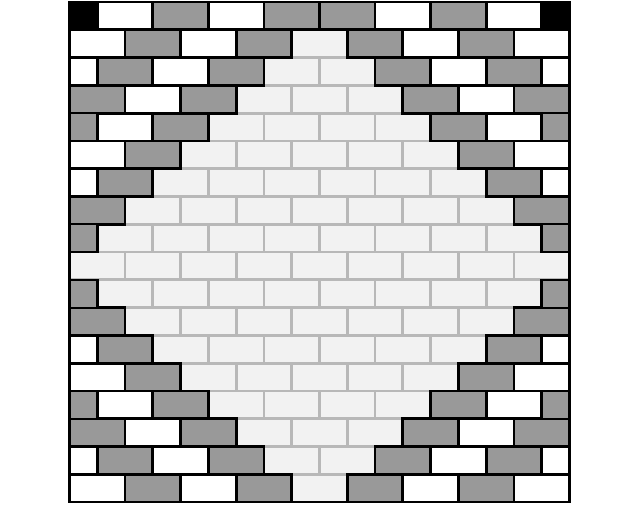%
      }
    \label{fig:evenEmptySet}}%
  \caption{Allowable diagonals shown in alternating grey and white,
    \emph{(a)}, for $\Tn{n}(l_i)$, where $(n,i)=(18,5)$, and
    \emph{(b)}, for $\Tn{18}(\varnothing)$.}
\label{fig:evenTn}
\end{figure}

We have $d_n(\lup{i})=n-i-1$, so we are interested in the number of
combinations of the above independently flippable diagonals with exactly
$k-(n-i-1)$ vertical dominoes.  That number is
\begin{align*}
  % Proofed June 13, 2012
  \ssum{k_1+k_2 = \\k -(n-i-1)} S(n-i-2,k_1)S(i-1,k_2).
\end{align*}
% Proofed June 13, 2012
The indices of the diagonals $\lup{i}$ for which $d_n(\lup{i}) \ge
d_n(\ldown{i})$ and $\rdown{i}$ for which $d_n(\rdown{i}) \ge
d_n(\rup{i})$, range from $1$ to $\flr{n-1}{2}$, as required for
(\ref{eq:tnka}).

% Table Proofed and Verified June 13, 2012 (now twice verified)

% VERIFIED FOR CORRECTNESS BY ALEJANDRO JAN 12, 2012
% Table Proofed and Verified June 13, 2012 (haha I hadn't see the above comment.  Oh well, twice verified is better than once :p)
\begin{table}[h]
  \label{tab:cases}
  \centering

  \begin{tabular}[h]{lllll}
    $\Tn{n}(a)$ &\multicolumn{4}{c}{Index and size of largest diagonal in this corner}\\
    \cline{2-5}
    $n$ even & $\ldown{j}$ ($j$ odd)& $\lup{j}$ ($j$ odd)& $\rdown{j}$ ($j$ even) & $\rup{j}$ ($j$ even)\\
    \hline
    $\Tn{n}(\lup{i})$& $j < i$ & $j>i$ & $j >i$ & c. 2(a)\\
    index $j$:& $ i-2$ & $i+2$ &$i+1$ & $i-1$\\
    \rowcolor{lightgrey}  size: & $i-2$ & $n-i-3$ & $n-i-2$ & $i-1$\\
    $\Tn{n}(\rdown{i})$ & c. 2(a) &$i<j$ & $j>i$ & $j<i$\\
    index $j$:& $i-1$ & $i+1$&$i+2$ & $i-2$ \\
    \rowcolor{lightgrey} size: &  $i-1$ & $n-i-2$&$n-i-3$& $i-2$\\
    $\Tn{n}(\ldown{i})$ & \multicolumn{4}{c}{Symmetric with $\Tn{n}(\rdown{i})$.}\\
    $\Tn{n}(\rup{i})$ & \multicolumn{4}{c}{Symmetric with $\Tn{n}(\lup{i})$.}\\
    \hline
    \hline
    $n$ odd & $\tleft{j}$ ($j$ even) & $\tright{j}$ ($j$ even)& $\bleft{j}$ ($j$ odd)& $\bright{j}$ ($j$ odd) \\
    \hline
    $\Tn{n}(\tright{i})$& $j < i$ & $j>i$ & $j < n-i-1$ & c. 2(a)\\
    index $j$:& $ i-2$ & $i+2$ &$n-i-2$ & $n-i$\\

    \rowcolor{lightgrey}  size: & $i-2$ & $n-i-3$ & $n-i-2$ & $i-1$\\
    $\Tn{n}(\bright{i})$ & $j < n-i-1$ &c. 2(a) & $j<i$ & $j>i$\\
    index $j$:&  $n-i-2$&$n-i$ &$i-2$ & $i+2$ \\
    \rowcolor{lightgrey}    size: &  $n-i-2$&$i-1$ & $i-2$& $n-i-3$\\
    $\Tn{n}(\tleft{i})$ & \multicolumn{4}{c}{Symmetric with $\Tn{n}(\tright{i})$.}\\
    $\Tn{n}(\bleft{i})$ & \multicolumn{4}{c}{Symmetric with $\Tn{n}(\bright{i})$.}\\
    \hline
  \end{tabular}

  \caption{The longest allowable diagonals in each of four corners for each $\Tn{n}(a)$.  Entries are calculated using the parity of $i$ and $j$, the avoidance of conflicts, and the requirement that $a$ be the longest diagonal in $\Tn{n}(a)$. Note that ``c. 2(a)'', above, refers to conflict 2(a) which occurs between diagonals $a$ and $b$ if $d_n(a) + d_n(b) \ge n$.}
  \label{tab:tnnnCases}
\end{table}

Now suppose $a=\varnothing$. If $i$ is the largest index such that
$d_n(\ldown{i}) < d_n(\lup{i})$ and $j$ is the largest index such that
$d_n(\rup{j}) < d_n(\rdown{j})$, then $\max(i,j) = \flr{n-2}{2}$ and
$|i-j|=1$. The allowable diagonals in $\Tn{n}(\varnothing)$ and their
respective sizes are shown in the table below (see \fref{evenEmptySet}).
\begin{center}
% VERIFIED FOR CORRECTNESS AND CONSISTENCY WITH TABLE 2 BY ALEJANDRO ON JAN 12, 2012
  \begin{tabular}{ll}
    % Proofed June 13, 2012
  $    \ldown{1},\ldown{3}, \ldots, \ldown{i} $&$  1,3, \ldots, i$\\
  $    \lup{i+2},\lup{i+4}, \ldots, \lup{n-3}$&$  n-i-3, n-i-5, \ldots, 2$\\
  $    \rup{2}, \rup{4}, \ldots,\rup{j}$&$  2, 4, \ldots,j,$\\
  $    \rdown{j+2}, \rdown{j+4}, \ldots,\rdown{n-2}$&$ n-j-3, n-j-5,
  \ldots,1.$
\end{tabular}
\end{center}

Choosing subsets of the independently flippable diagonals with $k$
vertical dominoes contributes the term
\begin{align*}
  % Proofed June 13, 2012
  \sum_{k_1+k_2 = k} S\left(\flr{n-2}{2},k_1\right)S\left(n- \left(\flr{n-2}{2}-1\right)-3,k_2\right),
\end{align*}
and since $n- \left(\lfloor (n-2)/2\rfloor-1\right)-3 = \lfloor (n-2)/2\rfloor$, this is
equal to (\ref{eq:tnkb}) for even $n$.

\textbf{Case $n$ odd:} The trivial covering is a vertical bond with
$\lfloor (n-2)/2\rfloor$ monominoes at the top (besides the two that
are fixed) and $\lceil (n-2)/2\rceil$ non-fixed monominoes along the
bottom boundary.  When diagonal $a$ is flipped, $d_n(a)$ horizontal
dominoes are added to the covering, instead of vertical dominoes.
Hence we argue for $H(n,k)$ rather than $V(n,k)$.

Now $\tleft{j}$ and $\tright{j}$ have even length, and $\bleft{j}$ and
$\bright{j}$ have odd length (see \tabref{tnnnCases}).  For example,
the allowable diagonals in $\Tn{n}(\tright{i})$ are shown in
\fref{oddExample} (for $(n,i)=(17,6)$), and their respective lengths
are

% VERIFIED FOR CORRECTNESS AND CONSISTENCY WITH TABLE 2 BY ALEJANDRO ON JAN 12, 2012
% Proofed June 13, 2012
\begin{center}
\begin{tabular}{ll}
  \centering
  $    \tleft{1},\tleft{3}, \ldots, \tleft{i-2} $&$  1,3, \ldots, i-2,$\\
  $    \tright{i+2},\tright{i+4}, \ldots, \tright{n-3}$&$  n-i-3, n-i-5, \ldots, 2,$\\
  $    \bleft{1}, \bleft{3}, \ldots,\bleft{n-i-2}$&$ 1,3,\ldots,n-i-2,$\\
  $    \bright{n-i}, \bright{n-i+2}, \ldots,\bright{n-2}$&$  i-1, i-3, \ldots,1.$
\end{tabular}
\end{center}

\begin{figure}[h]
  \centering
  \subfigure[]{%
    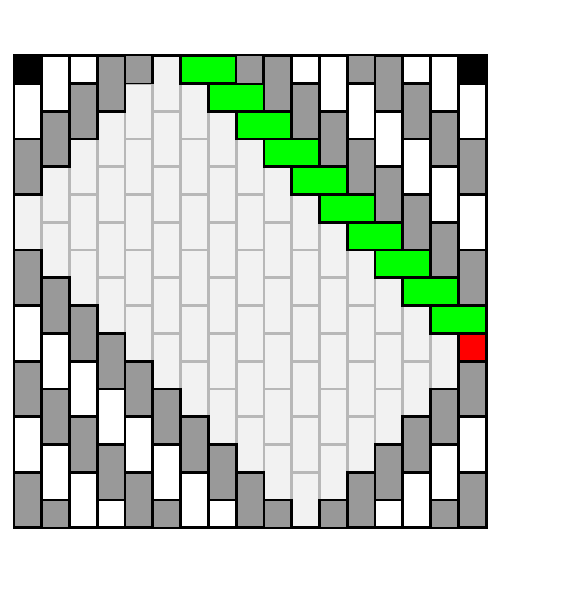%
    \label{fig:oddExample}} %
  \hspace{0.3in}%
  \subfigure[]{%
    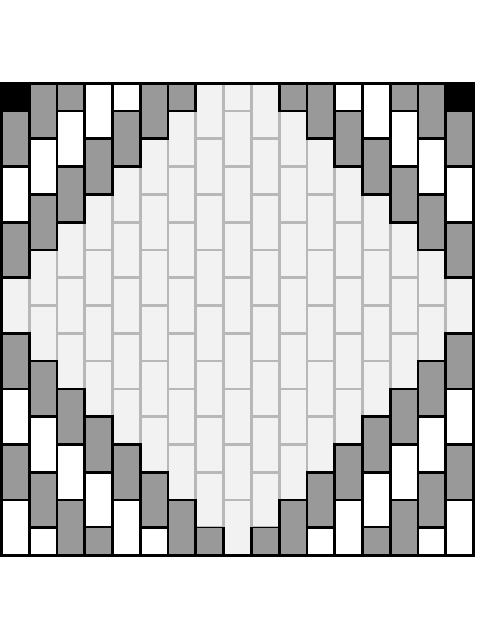%
    \label{fig:oddEmptySet}}%
  \caption{Allowable diagonals shown in alternating grey and white,
    \emph{(a)}, for $\Tn{n}(t_i)$, where $(n,i)=(17,6)$, and
    \emph{(b)}, for $\Tn{17}(\varnothing)$.}
\label{fig:oddTn}
\end{figure}

Once again $d_n(\tright{i})=n-i-1$, so we are interested in the number
of combinations of the above independently flippable diagonals with
exactly $k-(n-i-1)$ horizontal dominoes.  As before, that number is
\begin{align*}
  \ssum{k_1+k_2 = \\k -(n-i-1)} S(n-i-2,k_1)S(i-1,k_2).
\end{align*}

Now suppose $a=\varnothing$, then if $i$ is the largest index such
that $d_n(\tright{i}) < d_n(\tleft{i})$ and $j$ is the largest index
such that $d_n(\bright{j}) < d_n(\bleft{j})$ then $\max(i,j) =
\flr{n-2}{2}$ and $|i-j|=1$. The allowable leftward diagonals in
$\Tn{n}(\varnothing)$ and their respective sizes are given in the
table below.
 \begin{center}
\begin{tabular}{ll}
  \centering
  % Proofed June 13, 2012
  $    \tleft{2}, \tleft{4}, \ldots,\tleft{j}$&$  2, 4, \ldots,j,$\\
  $    \bleft{1},\bleft{3}, \ldots, \bleft{i} $&$  1,3, \ldots, i$\\
\end{tabular}
\end{center}
and by horizontal symmetry, the rightward diagonals have the same
lengths.  We conclude that the coverings with $k$ horizontal dominoes of
$\Tn{n}(\varnothing)$ is also generated by (\ref{eq:tnkb}) when $n$ is
odd.\qed

The terms $\vh(n,k)z^k$ can be summed over $k$ to obtain the
generating polynomial $T(n,z)$ (same as $\vh_n(z)$), mentioned in
Conjecture 4 of \cite{EricksonRuskeySchurch2011}.

\begin{theorem}
  \label{thm:vdominoes}
  The generating polynomial for (\ref{eq:tnka}-\ref{eq:tnkb}) is
  % VERIFIED AND CORRECTED APRIL 13
  \begin{align}
    \label{eq:tnz}
    \vh_n(z):=2\sum_{i=1}^{\flr{n-1}{2}}
    S_{n-i-2}(z)S_{i-1}(z)z^{n-i-1} +
    \left(S_{\flr{n-2}{2}}(z)\right)^2,
  \end{align}
  where $S_n(z) = \sum_{k\in \mathbb Z} S(n,k)z^k$.  This
  ``generates'' $V(n,k)$ for even $n$, and $H(n,k)$ for odd $n$.
\end{theorem}
%Upward diagonals have even length, and downward diagonals have odd
%length.
\proof This follows from \lref{vdominoes}.  \qed

The degree of $\vh_n(z)$ is $\frac{n^2-n}{2}-(n-1)$, because this is
the largest number of vertical dominoes possible in a covering of
$\Tn{n}$, for even $n$ (and horizontal dominoes for odd $n$).  For
example, the covering with all $l_i$ flipped up and all $r_i$ flipped
down contains exactly $n-1$ horizontal dominoes.

\newcommand{\coeff}[1]{\langle z^{#1}\rangle}

The coefficients of $\vh_n(z)$ are listed in \tabref{VH} up to $n =
10$, and the following conjecture is true at least up to $n=20$.  If
$Q(z)$ is a polynomial, then write $\coeff{k}Q(z)$ to denote the
coefficient of $z^k$.

\begin{conjecture}~ %non breaking space puts itemize on next line
  % Proofed June 13, 2012
  \label{conj:vhcon}
  \begin{enumerate}[(a)]
  \item For $k\le n-2$, we have $\coeff{k}\vh_n(z) =
    \coeff{k}\prod_{m\ge 0}(1+z^m)^2$, the number of partitions of $k$
    into distinct parts with two types of each part (see A022567 in
    \cite{Sloane}).
  \item For $0\le k < n-3$, we have
    \begin{align*}
      \coeff{\deg\left(\vh_n(z)\right)-k}\vh_n(z) =
      2\coeff{k}\prod_{m\ge 0}(1+z^m),
    \end{align*}
    twice the number of partitions of $k$ into
      distinct parts (see A000009 in \cite{Sloane}).
  \end{enumerate}
\end{conjecture}

% entries sight verified for n=10 and most of n=9, as well as
% n=2,3,4,5 and some of 6 against big table in previous version of
% paper. Sept 23, 2012
\begin{table}[h]\centering{\tiny \begin{tabular}[h]{l|rrrrrrrrrrrrrrrr}$n\backslash z^k$& 0& 1& 2& 3& 4& 5& 6& 7& 8& 9& 10& 11& 12& 13& 14& 15\\\hline2&   1&&&&&&&&&&&&&&&\\3&   1&   2&&&&&&&&&&&&&&\\4&   1&   2&   3&   2&&&&&&&&&&&&\\5&   1&   2&   3&   6&   4&   2&   2&&&&&&&&&\\6&   1&   2&   3&   6&   9&   8&   7&   6&   2&   2&   2&&&&&\\7&   1&   2&   3&   6&   9&  14&  15&  14&  14&  10&   8&   6&   4&   2&   2&   2\\8&   1&   2&   3&   6&   9&  14&  22&  24&  25&  28&  25&  22&  19&  14&  10&  10\\9&   1&   2&   3&   6&   9&  14&  22&  32&  37&  42&  49&  48&  49&  46&  38&  34\\10&   1&   2&   3&   6&   9&  14&  22&  32&  46&  56&  66&  78&  84&  90&  92&  88\\$n\backslash z^k$& 16& 17& 18& 19& 20& 21& 22& 23& 24& 25& 26& 27& 28& 29& 30& 31\\\hline8&   8&   4&   4&   2&   2&   2&&&&&&&&&&\\9&  30&  24&  20&  16&  12&  12&  10&   6&   4&   4&   2&   2&   2&&&\\10&  81&  76&  69&  58&  51&  44&  38&  34&  28&  22&  20&  16&  14&  12&   8&   6\\$n\backslash z^k$& 32& 33& 34& 35& 36& 37& 38& 39& 40& 41& 42& 43& 44& 45& 46& 47\\\hline10&   4&   4&   2&   2&   2&&&&&&&&&&&\\\end{tabular}}\caption{Table of coefficients of $\vh_n(z)$ for $2\le n \le 10$.  The $(n,k)$th entry represents the number of coverings of $\Tn{n}$ with $k$ vertical dominoes when $n$ is even, and $k$ horizontal dominoes when $n$ is odd.}\label{tab:VH}\end{table}

\begin{sloppypar}
  Rotating a covering of $\Tn{n}$ by $\pi/2$ radians interchanges
  vertical and horizontal dominoes, and this transformation can be
  applied to the generating polynomial $\vh_n(z)$ to obtain the
  polynomial $\vh_n(z^{-1})z^{\frac{n^2-n}{2}}$.  Thus we can easily
  derive the bivariate generating polynomial $R_n(x,y)$, whose
  coefficient of $x^vy^h$ is the number of tatami coverings with
  exactly $v$ vertical dominoes and $h$ horizontal dominoes.
\end{sloppypar}

Our remarks prove the following corollary.

\begin{corollary}
  Let $R_n(x,y)$ be as defined above.  We have
  \begin{align}
    \label{eq:rnz}
    R_n(x,y) = 2\vh_n(xy^{-1})y^{\frac{n^2-n}{2}} +
    2\vh_n(x^{-1}y)x^{\frac{n^2-n}{2}}.
  \end{align}
\end{corollary}

% not important enough to warrant a whole page
%\input{Rtable.tex}

We list some basic properties of $R_n(x,y)$.
\begin{itemize}
  % Proofed June 13, 2012
\item The degree of $R_n(x,1)$ as well as the degree of every term in
  $R_n(x,y)$ is $\frac{n^2-2}{2}$;
\item the polynomial $R_n(x,y)$ can be recovered from $R_n(x,1)$, and
  the latter is the generating polynomial for the set of all $n\times
  n$ coverings with $n$ monominoes with exactly $v$ vertical dominoes
  (or $h$ horizontal dominoes);
\item the polynomial $R_n(x,1)$ is self reciprocal because of
  interchangeability of vertical and horizontal dominoes; and finally,
\item the polynomial $R_n(x,1)$ has similar properties to those listed
  for $\vh_n(z)$ in \conref{vhcon}, in the sense that for some
  increasing integer function $f$, we have $\langle x^k\rangle
  R_{n}(x,1) = \langle x^k \rangle R_{n+1}(x,1)$, whenever $k<
  f(n+1)$.
\end{itemize}

% We observe again that, as $n$ grows, $R_n(x,1)$ appears to converge to
% a series, where $[k]R_n(x,1) = [k]R_N(x,1)$ for $k\le n-2$ and $N>n$.
% The series stabilizing in the reciprocal is the same, because
% $R_n(x,1)$ is self-reciprocal.  Reading from $R_{20}(z)$ in
% \tabref{Rn}, we have
% \begin{align*}
%   2,4,6,12,18,28,44,64,92,132,186,256,352,476,638,852,1124,
% \end{align*}
% and this is simply double the sequence observed in $\vd_n(z)$.

If there is an even number of dominoes, which is the case when
$(n^2-n)/4$ is an integer, then $\langle x^ky^k \rangle R_n(x,y) =
4\coeff{k}\vh_n(z)$, where $k=(n^2-n)/4$. Rotating the covering maps
$k$ vertical dominoes to $k$ horizontal dominoes, and vice versa.  The
coverings counted by these coefficients are called \emph{balanced
  tatami coverings}, appropriately named by Knuth (private
communication), because the number of vertical and horizontal dominoes
are equal.  Here is $\coeff{k}\vh_n(z)$ for $2\le n \le 56$: $0$, $
0$, $ 2$, $ 2$, $ 0$, $ 0$, $ 10$, $ 20$, $ 0$, $ 0$, $ 114$, $ 210$,
$ 0$, $ 0$, $ 1322$, $ 2460$, $ 0$, $ 0$, $ 16428$, $ 31122$, $ 0$, $
0$, $ 214660$, $ 410378$, $ 0$, $ 0$, $ 2897424$, $ 5575682$, $ 0$, $
0$, $ 40046134$, $ 77445152$, $ 0$, $ 0$, $ 563527294$, $ 1093987598$,
$ 0$, $ 0$, $ 8042361426$, $ 15660579168$, $ 0$, $ 0$, $
116083167058$, $ 226608224226$, $ 0$, $ 0$, $ 1691193906828$, $
3308255447206$, $ 0$, $ 0$, $ 24830916046462$, $ 48658330768786$, $
0$, $ 0$, $366990100477712$, (see A182107 in \cite{Sloane}). Note that
this is perhaps better viewed as four sequences, one for each ${0\le j
< 4}$ such that $n \pmod 4 = j$.

% \Begin{table}[h]
%   \centering
%   {\small
% \begin{tabular}{l|llllllllllll}
%   n&  2 &3 &4 &5 &6 &7 &8 &9 &10 &11 &12 &13 \\
%   $b_n$&  0 &0 &2 &2 &0 &0 &10 &20 &0 &0 &114 &210\\
%   \hline
% n&  14 &15 &16 &17 &18 &19 &20 &21 &22 &23 &24\\
% $b_n$&  0 &0 &1322 &2460 &0 &0 &16428 &31122 &0 &0 &214660\\
%   \hline
% n&  25 &26 &27 &28 &29 &30 &31 &32 &33 &34 &35\\
% $b_n$&  410378 &0 &0 &2897424 &5575682 &0 &0 &40046134 &77445152 &0 &0 \\
% %  \vspace{0.2cm}
% \end{tabular}
% }
% \caption{The sequence $b_n$ gives the number of balanced coverings for each $2\le n \le 61$, up to rotational symmetry.  Thus $b_n = [z^k]\vh_n(z)$, where $k=(n^2-n)/4$, if $k$ is an integer, and otherwise $b_n=0$.}
% \label{tab:balanced}
% \end{table}

\section{Exhaustive generation of tatami coverings in constant amortized time}
\label{sec:combalg}
\newcommand{\cc}{\texttt{C4}}%
\newcommand{\modcc}{\texttt{modC}}%
\newcommand{\Sset}{\mathbf S}

In this section, we present procedure \texttt{genVH}$(n,k)$, which
generates the coverings counted by $\vh(n,k)$, while doing a constant
amount of data structure change per covering that is produced (example
output is shown in \fref{n8k7}). Let $\Sset(n,k)$ denote the set of
subsets of $\set{1,2,\ldots, n}$ whose elements sum to $k$; thus
$|\Sset(n,k)| = S(n,k)$.  The procedure follows naturally from the
sums in Equation (\ref{eq:tnka}-\ref{eq:tnkb}), since each term
$S(a,i)S(b,j)$ counts some set of \fd~and \bd-oriented diagonals. The
sets $\Sset(a,i) \times \Sset(b,j)$ are generated in constant
amortized time (CAT) by a modification of \cc~from
\cite{BaronaigienRuskey1993} (see \lstref{modc}). Our modified
algorithm, \modcc, is invoked for each sum term of Equation
(\ref{eq:tnka}-\ref{eq:tnkb}).  Procedure \modcc~is CAT for the same
reasons that \cc~is CAT.

\begin{figure}[h]
  \centering
\includegraphics[width=\textwidth]{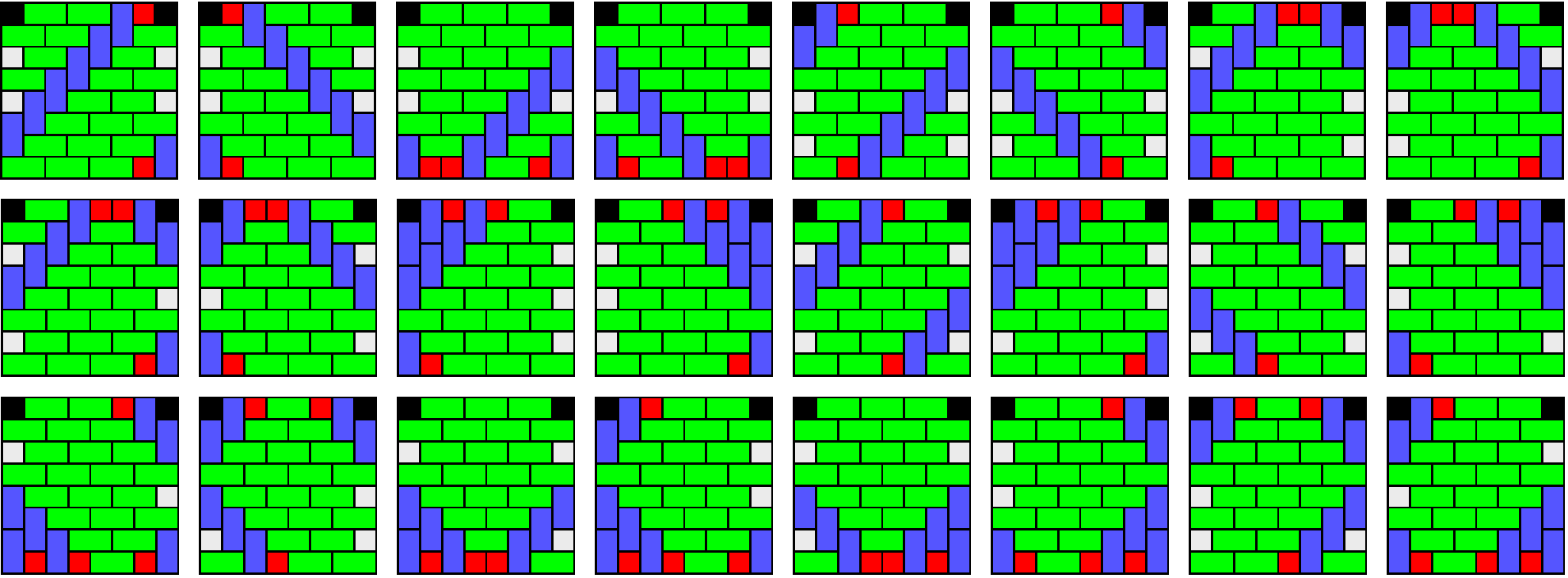}
\caption{The coverings of $\Tn{8}$ with exactly 7 vertical dominoes.  This
  is the output of \texttt{genVH}$(8,7)$ printed in the order the
  coverings are generated (as one would naturally read text).}
\label{fig:n8k7}
\end{figure}

\begin{lstlisting}[float,caption = {Python code for a modified version
of \cc~from \cite{BaronaigienRuskey1993} to compute $\Sset(a,i)\times
\Sset(b,j)$.  Global variables \texttt{aiSet} and \texttt{bjSet} are
the lists representing $\Sset(a,i)$ and $\Sset(b,i)$, respectively},
label={lst:modc}]
def modC(a,i,b,j,comp,isFirst):
    global aiSet,bjSet
    if( a == 0):
        if(isFirst):
            modC(b,j,0,0,False,False)
        else:
            Output(aiSet,bjSet)
    else:
        if(isFirst):
            L = aiSet
        else:
            L = bjSet
        if( i > a*(a+1)/2 ):
            i = a*(a+1)/2 - i; comp = not comp
        if( i<a ):
            if(comp):
                L[a] = L[0]; L[0] = i+1
                modC( i, i, b, j, comp, isFirst)
                L[0] = L[a]; L[a] = a+1
            else:
                modC( i, i, b, j, comp, isFirst)
        else:
            L[a] = L[0]; L[0] = a
            if(comp):
                modC( a-1, i, b, j, comp, isFirst)
                L[0] = L[a]; L[a] = a+1
                modC( a-1, i-a, b, j, comp, isFirst)
            else:
                modC( a-1, i-a, b, j, comp, isFirst)
                L[0] = L[a]; L[a] = a+1
                modC( a-1, i, b, j, comp, isFirst)
\end{lstlisting}

There is one subtlety involved in exploiting the CATness of \cc.
Invoking \cc$(a,i)$ requires $\Omega(a)$ preprocessing steps if its
input list is recreated for each call, but \cc$(a,i)$ may not produce
so many combinations for small $a$ and large $i$.  The result is that
we may make many calls to \modcc~ that require too much preprocessing,
but this is dealt with, as follows: a top level call to \cc$(i,j)$ in
\cite{BaronaigienRuskey1993} takes the list $[i+1,2,3,\ldots, i+1]$,
which requires $i+1$ steps to create, however, \cc$(i,j)$ also
concludes with the same list (see \lstref{modc}).  Let $A$ and $B$ be
the largest integers for which \modcc~is called to compute $\Sset(a,i)
\times \Sset(b,i)$.  We set~\texttt{aiSet}$ = [1,2,3,\ldots, A]$ and
\texttt{biSet}$ = [1,2,3,\ldots, B]$, and by
setting~\texttt{aiSet}$[0]=a+1$ and~\texttt{biSet}$[0]=b+1$, we
initialize for each call to \modcc~with exactly two operations.

\begin{theorem}
  \label{thm:cattnk}
  The coverings in $\Tn{n}$ with exactly $k$ vertical dominoes if $n$ is
  even and horizontal dominoes if $n$ is odd, can be exhaustively
  generated in constant amortized time.
\end{theorem}
\proof The outer procedure does a constant amount of work per call to
\modcc.  This subroutine is CAT, so the outer procedure is also
CAT.\qed

\section{A mysterious factor of $\vh_n(z)$}
\label{sec:pn}

In this section we prove that the generating polynomial $\vh_n(z)$ has
(very nearly) the factorisation conjectured in
\cite{EricksonRuskeySchurch2011}.  We use the following lemma.

\begin{lemma}
  % not in GKP
  \label{lem:flrxapp} For all $x\ge 0$,
  \begin{align}
  \label{eq:flrxapp}
    \lfloor x \rfloor = \sum_{k\ge 1} \left\lfloor \frac{x}{2^k}+\frac{1}{2} \right\rfloor.
  \end{align}
\end{lemma}
\proof Let $n = \lfloor x \rfloor$ and apply strong induction on $n$.
Clearly \eref{flrxapp} holds for the base case, when $n=0$.  Suppose
it holds for $0,1,\ldots, n-1$, then
\begin{align*}
  % Proofed/Verified June 13, 2012
  \sum_{k\ge 1} \round{x}{2^k} =& \round{x}{2} + \sum_{k\ge 1}
  \round{\frac{x}{2}}{2^k} \\
  =& \round{x}{2} + \flr{x}{2} = \round{\lfloor x \rfloor}{2} +
  \flr{\lfloor x \rfloor}{2} = \lfloor x \rfloor.
\end{align*}
%[BELOW] VERIFIED BY ALEJANDRO ON JULY 7, 2012
Two applications of Equation (3.11) in \cite{GrahamKnuthPatashnik1994}
yield the penultimate equation and the final equation follows by
considering the parity of $\lfloor x \rfloor$, or by using (3.26) in
\cite{GrahamKnuthPatashnik1994} with $m=2$ and $x'=x/2$.\qed

%%%%%%%%%%%%%%%%%%%%%%%%%%%%%%%%
%CALCULATION SHOULD BE HAND VERIFIED ONE MORE TIME MARCH 30
\begin{theorem}
  \label{thm:cyclotomic}
  The generating polynomial $\vh_n(z)$ has the factorisation
  \begin{align*}
    \vh_n(z) = P_n(z)D_n(z)
  \end{align*}
  where $P_n(z)$ is a polynomial and
  \begin{align}
    \label{eq:dn}
  D_n(z) = \prod_{j\ge 1} S_{\left\lfloor \frac{n-2}{2^j}
      \right\rfloor}(z).
  \end{align}
\end{theorem}
\proof % Proofed June 13, 2012
We prove that $D_n(z)$ divides $\vh_n(z)$ by using the
factorisation of $S_n(z)$ into cyclotomic polynomials
(\cite{EricksonRuskeySchurch2011}, Lemma 5),
\begin{align}
  \label{eq:snphi}
  S_n(z) = \prod_{j\ge 1} \Phi_{2j} (z)^{\lfloor \frac{n+j}{2j}
    \rfloor}.
\end{align}
and showing that the power of $\Phi_i(z)$ is greater in each term of
$\vh_n(z)$ than it is in $D_n(z)$.

The power of $\Phi_{2j}(z)$ in $D_n(z)$ is obtained by substituting
\eref{snphi} into \eref{dn},
\begin{align*}
  D_n(z)= \prod_{i\ge 1}S_{\flr{n-2}{2^i}}(z) =& \prod_{i\ge 1,j\ge 1}
   \Phi_{2j} (z)^{\left\lfloor
      \frac{\flr{n-2}{2^i}+j}{2j}
    \right\rfloor} \\
  =& \prod_{j\ge 1}  \Phi_{2j} (z)^{\sum_{i\ge 1}
    \flr{\flr{n-2}{2^i}+j}{2j}}.
\end{align*}

We simplify $D_n(z)$ to
  \begin{align}
    \label{eq:dnFactors}
  D_n(z)=& \prod_{j\ge 1}  \Phi_{2j} (z)^{\flr{n-2}{2j}},
\end{align}
by applying \lref{flrxapp} and with Equation (3.11) in
\cite{GrahamKnuthPatashnik1994}.

Expanding the second term of $\vh_n(z)$ gives
\begin{align*}
  \left(S_{\frac{n-2}{2}}(z)\right)^2 = \prod_{j\ge 1} 
    \Phi_{2j} (z)^{2\flr{\frac{n-2}{2}+j}{2j}},
\end{align*}
which is divisible by $D_n(z)$, since
\begin{align*}
  \flr{n-2}{2j}\le 2\flr{\frac{n-2}{2}+j}{2j},
\end{align*}
for all $j\ge 1$ and positive even integers $n$.

The other terms in $\vd_n(z)$ are of the form
\begin{align*}
  S_{n-k-2}(z)S_{k-1}(z)z^d = \left(\prod_{j>0}
    \Phi_{2j}(z)^{\flr{(n-k-2) + j}{2j}}\right)\left(\prod_{j>0}
    \Phi_{2j}(z)^{\flr{(k-1)+j}{2j}}\right)
  z^d
\end{align*}
for each $1\le k\le \flr{n-1}{2}$
where $d$ is the appropriate power of $z$. These terms are all
divisible by $D_n(z)$ if the exponents in \eref{dnFactors} satisfy
\begin{align}
  \label{eq:firstTermIneq}
  \flr{n-2}{2j}\le
  \left\lfloor \frac{k-1}{2j} + \frac{1}{2} \right\rfloor
  +\left\lfloor \frac{n-k-2}{2j} + \frac{1}{2} \right\rfloor.
\end{align}
Let $r_1$ and $r_2$ be integers such that ${0\le r_i < 2j}$ and
$\frac{k-1}{2j} = \flr{k-1}{2j} + \frac{r_1}{2j}$ and $\frac{n-2}{2j}
= \flr{n-2}{2j} + \frac{r_2}{2j}$.  We eliminate occurrences of
$\flr{k-1}{2j}$ and $\flr{n-2}{2j}$ from Inequality
(\ref{eq:firstTermIneq}), since they are integers and can be removed
from floors, and rewrite the inequality as
\begin{align}
  \label{eq:firstTermIneqRewrite}
  0 \le \left\lfloor \frac{r_1}{2j} + \frac{1}{2} \right\rfloor
  +\left\lfloor \frac{r_2-r_1-1}{2j} + \frac{1}{2} \right\rfloor.
\end{align}
It is straightforward to show that if the second term is $-1$, then
the first term is equal to $1$.

Therefore, $D_n(z)$ divides each and every term of $\vh_n(z)$.\qed
% end of theorem proof

% entries sight verified for most of n=3,4,5,11 and a few spot checks
% against big table in previous version of paper. Sept 23, 2012
\begin{table}[h]\centering{\tiny \begin{tabular}[h]{l|rrrrrrrrrrrrrrrr}$n\backslash z^k$& 0& 1& 2& 3& 4& 5& 6& 7& 8& 9& 10& 11& 12& 13& 14& 15\\\hline3&   1&   2&&&&&&&&&&&&&&\\4&   1&   1&   2&&&&&&&&&&&&&\\5&   1&   1&   2&   4&   0&   2&&&&&&&&&&\\6&   1&   0&   1&   2&   2&  -2&   2&&&&&&&&&\\7&   1&   0&   1&   2&   2&   4&  -2&   4&   0&   2&  -2&   2&&&&\\8&   1&   0&   1&   1&   2&   3&   4&  -2&   2&   0&   4&  -2&   2&  -2&   2&\\9&   1&   0&   1&   1&   2&   3&   4&   6&  -2&   6&   0&   8&  -2&   4&  -4&   6\\10&   1&  -1&   1&   0&   1&   1&   1&   2&   2&  -6&   6&  -2&   6&  -6&   4&  -4\\11&   1&  -1&   1&   0&   1&   1&   1&   2&   2&   4&  -8&  10&  -4&  10&  -8&   8\\$n\backslash z^k$& 16& 17& 18& 19& 20& 21& 22& 23& 24& 25& 26& 27& 28& 29& 30& 31\\\hline9&  -2&   4&  -2&   2&  -2&   2&&&&&&&&&&\\10&   6&  -6&   6&  -4&   4&  -4&   2&&&&&&&&&\\11&  -8&  10& -10&  12&  -8&  10& -12&  10&  -6&   6&  -6&   6&  -4&   4&  -4&   2\\\end{tabular}}\caption{Table of coefficients of $P_n(z)$ for $3\le n \le 11$.}\label{tab:P}\end{table}

Our computer investigations show that $P_n(z)$ is irreducible for ${1<n<200}$,
and we know the complete factorisation of $S_k(z)$, for each positive
integer $k$.  We suspect, therefore, that the complete factorisation is
\begin{align}
  \label{eq:completeFactorisation}
  \vh_n(z) = P_n(z)\prod_{j\ge 1} \Phi_{2j} (z)^{\flr{n-2}{2j}}.
\end{align}

The factor $P_n(z)$ is somewhat more mysterious than $D_n(z)$; e.g.,
we have no formula to express it besides $\vh_n(z)/D_n(z)$. Take
$P_{11}(z)$ for example, which is equal to
% {\small
% \begin{tabular}{llllllllllllllll}
%   1&  -1z^{1}&   1z^{2}&   0z^{3}&   1z^{4}&   1z^{5}&   1z^{6}&   2z^{7}&   2z^{8}&   4z^{9}&  -8z^{10}&  10z^{11}&  -4z^{12}&  10z^{13}&  -8z^{14}&  8z^{15}\\
%   -8z^{16}&  10z^{17}& -10z^{18}& 12z^{19}& -8z^{20}&  10z^{21}& -12z^{22}&  10z^{23}&  -6z^{24}&   6z^{25}&  -6z^{26}&   6z^{27}&  -4z^{28}&   4z^{29}&  -4z^{30}&   2z^{31}.&
% \end{tabular}
% }
$ 1-1z^{1}+ 1z^{2}+ 0z^{3}+ 1z^{4}+ 1z^{5}+ 1z^{6}+ 2z^{7}+ 2z^{8}+
4z^{9}-8z^{10}+ 10z^{11}-4z^{12}+ 10z^{13}-8z^{14}+ 8z^{15} -8z^{16}+
10z^{17}-10z^{18}+ 12z^{19}-8z^{20}+ 10z^{21}-12z^{22}+
10z^{23}-6z^{24}+ 6z^{25}-6z^{26}+ 6z^{27}-4z^{28}+ 4z^{29}-4z^{30}+
2z^{31}$. The coefficients are almost all non-zero, the central
coefficients are larger than the ones at the tails, a great many of
them are even, they alternate in sign for a long stretch, and the
polynomial is irreducible.

The degree of $P_{11}(z)$ is $\deg(P_{11}(z)) =\deg(\vh_{11}(z))-
\deg(D_{11}(z))$, both of which are easily calculated.  In general
$\deg(P_n(z))$ is equal to the sum of the sequence of largest odd
divisors of the numbers $1,2, \ldots, n-2$, which is a sequence with
some nice properties (see A135013 in \cite{Sloane}).

\begin{theorem}
  \label{thm:degpn}
   For each $n\ge 2$,
\begin{align*}
  \deg(P_n(z)) = \sum_{k=1}^{n-2} \op(k),
\end{align*}
where $\op(k)$ is the largest odd divisor of $k$.
\end{theorem}

\proof   % [WHOLE PROOF VERIFIED BY ALEJANDRO APRIL 18, 2012]
% and again June 13, 2012
\tref{cyclotomic} gives the degree of $D_n(z)$ so we can write
\begin{align}
  \label{eq:degp}
  \deg(P_n(z)) = \binom{n-1}{2} - \sum_{k\ge 1}
  \binom{\flr{n-2}{2^k}+1}{2},
\end{align}
since $\deg(S_n(z)) = \binom{n+1}{2}$.

The proof that $\sum_{k=1}^{n} \op(k) = \deg(P_{n+2}(z))$ is by induction, and the
base case, where $n=0$, is easily verified.  Let $p_{n} =
\deg(P_{n}(z))$, for $n\ge 2$, to abbreviate the notation.  It remains
for us to show that $p_{n+3} - p_{n+2} = \op(n+1)$.

Let $n'01^\alpha$ be the binary representation of $n$, so that $n+1 =
n'10^\alpha$, and let $\I{A}=1$ if the statement $A$ is true, and
$\I{A}=0$ otherwise.
% $n=0^{\alpha-k}1n' = 1 + 1^{\alpha-k}0n' 1 + \flr{n}{2^k} , &\text{
%   if } k \le \alpha$
%
% $\flr{n'}{2^{k-\alpha-1}}= \flr{n}{2^k} , &\text{ if } k > \alpha. $
Observe that
\begin{align}
%[VERIFIED BY ALEJANDRO APRIL 18]
  \label{eq:flr2k}
  \flr{n+1}{2^k} = \I{k\le \alpha} + \flr{n}{2^k},
\end{align}
 which we use to simplify
\begin{align*}
  \sum_{k \ge 1} \left(\binom{\flr{n+1}{2^k}+1}{2} -
    \binom{\flr{n}{2^k}+1}{2} \right),
\end{align*}
and write
\begin{align*}
  % [VERIFIED BY ALEJANDRO APRIL 18]
  p_{n+3} - p_{n+2} =  (n+1) - \sum_{k=1}^{\alpha}\left( \flr{n}{2^k} +1 \right).
\end{align*}
Using \eref{flr2k} and the fact that $(n+1)/2^k$ is an integer for $1\le
k\le \alpha$, we write
\begin{align*}
  % [VERIFIED BY ALEJANDRO APRIL 18]
  p_{n+3} - p_{n+2} =& (n+1) - \sum_{k=1}^{\alpha}\left( \round{n+1}{2^k} \right),
\end{align*}
and then express this as the remaining sum terms in \eref{flrxapp}
\begin{align*}
  p_{n+3} - p_{n+2} =& \sum_{k\ge \alpha + 1} \left(  \round{n+1}{2^k} \right)\\
  =& \sum_{k-\alpha \ge 1} \left( \left\lfloor
      \frac{\frac{n+1}{2^\alpha}}{2^{k-\alpha}}+ \frac{1}{2} \right\rfloor \right).
\end{align*}

Applying \eref{flrxapp} again, we have $ p_{n+3} - p_{n+2} =
(n+1)/2^\alpha$, which is equal to $\op(n+1)$, as required.\qed

In addition to finding $\deg(P_{11}(z))$, we can evaluate at $z=1$
with $P_{11}(1) = \vh_{11}(1)/D_{11}(1) = 22$, a ratio which is also
easy to calculate in general because $\vh_{n}(1)$ and $D_n(1)$ have
well understood combinatorial interpretations. It also leads to an
interesting sequence, whose derivation for all $n$ is given below.

\begin{theorem}
  \label{thm:nu}
  The sum of the coefficients of $P_n(z)$ is equal to
  $n2^{\nu(n-2)-1}$, where $\nu(n)$ is the number of $1$s in the
  binary representation of $n$
\end{theorem}
\proof The sum of the coefficients of $P_n(z)$ is equal to $P_n(1)$,
which is expressible as $\vh_{n}(1)/D_{n}(1)$.  The numerator
evaluates to $n2^{n-3}$, since this is the number of coverings in
$\Tn{n}$, and the denominator is evaluated as described below.

It is well known that $\Phi_{k}(1) = p$ if $k$ is a non-zero power of
a prime $p$ and $\Phi_k(1) = 1$ if $k$ is divisible by two distinct
primes (see \cite{Lang1970}, p.74). We can evaluate $D_n(1)$ using
\eref{dnFactors},
\begin{align*}
  D_n(1)= \prod_{i\ge 1} \Phi_{2i}(1)^{\flr{n-2}{2i}} =  2^{\sum_{i\ge 1}\flr{n-2}{2^i}},
\end{align*}
by ignoring the factors for which $2i$ is not a power of $2$.  Apply
Equation (4.24) in \cite{GrahamKnuthPatashnik1994} to obtain $D_n(1) =
2^{n-2 - \nu(n-2)}$. Thus
\begin{align*}
  P_n(1) = n2^{n-3 - (n-2) + \nu(n-2)} = n 2^{\nu(n-2)-1}.
\end{align*}\qed

We have verified that $P_n(z)$ is irreducible over the integers for
${1<n < 200}$, but we do not understand its structure well enough to
prove it for all $n$.  We state below some of the observable structure
which has also been verified for ${1 < n < 200}$, as \conref{Pcon}, and
we plot some complex roots for odd $n$ up to $67$ in
\fref{complexPlotP}.

\begin{conjecture}
  \label{conj:Pcon}
~
  \begin{enumerate}[(a)]
    % this one is false
    % \item The sequence $\set{[z^i]P_n(z)}_{i\ge n-1}$ alternates
    %   between non-positive and non-negative integers, beginning with
    %   the former.
  \item If $k \ge 1$ and $n \pmod{2^k} = 2$, then $\coeff{i}P_n(z) = \coeff{i}P_{n+j}(z)$
    for $i\le \frac{n-2}{2^{k-1}}$ and $j\le 2^k$.
  \item When $n$ is odd, $P_n(z)$ has exactly one real root
    $\alpha_n$, with ${-1< \alpha_n \le -0.5}$, and $\set{\alpha_n}_{n
      \text{ odd }}$ is a monotonically decreasing sequence.
  \item When $n$ is even, $P_n(z)$ has no real root.
    % \item The sequences of polynomials $\set{P_n(z)}_{n \text{ even
    %     }}$ and $\set{P_n(z)}_{n \text{ odd }}$ are also convergent
    %   [SEE FIGURE BLAH].  [Probably not true, judging from plots]
  \item The polynomial $P_n(z)$ is irreducible over the integers for
    $n\ge 2$.
    % TOO OBVIOUS IN LIGHT OF THE PLOTS.
  % \item The polynomial $P_n(z)$ has exactly one local minimum, one
  %   local maximum, and one inflection point for $n\ge 7$, when $n$ is
  %   odd.  When $n$ is even it has only the first of these three, for
  %   $n\ge 4$.
  \item The alternating sums of coefficients are given by the
    generating function
    \begin{align}
      \label{eq:genpn1}
      \sum_{n\ge 2} P_n(-1)z^{n-2} =
      \frac{(1+z)(1-2z)}{(1-2z^2)\sqrt{1-4z^2}}.
    \end{align}
  \item \label{enum:alternating} For even $n$, the sum of the absolute
    values of coefficients of $P_n(z)$ is equal to $P_n(-1)$ when
    $n\ge 20$.
  \end{enumerate}
\end{conjecture}

\begin{figure}[ht]
%  \begin{threeparttable}
\centering
  \includegraphics[width=0.8\textwidth]{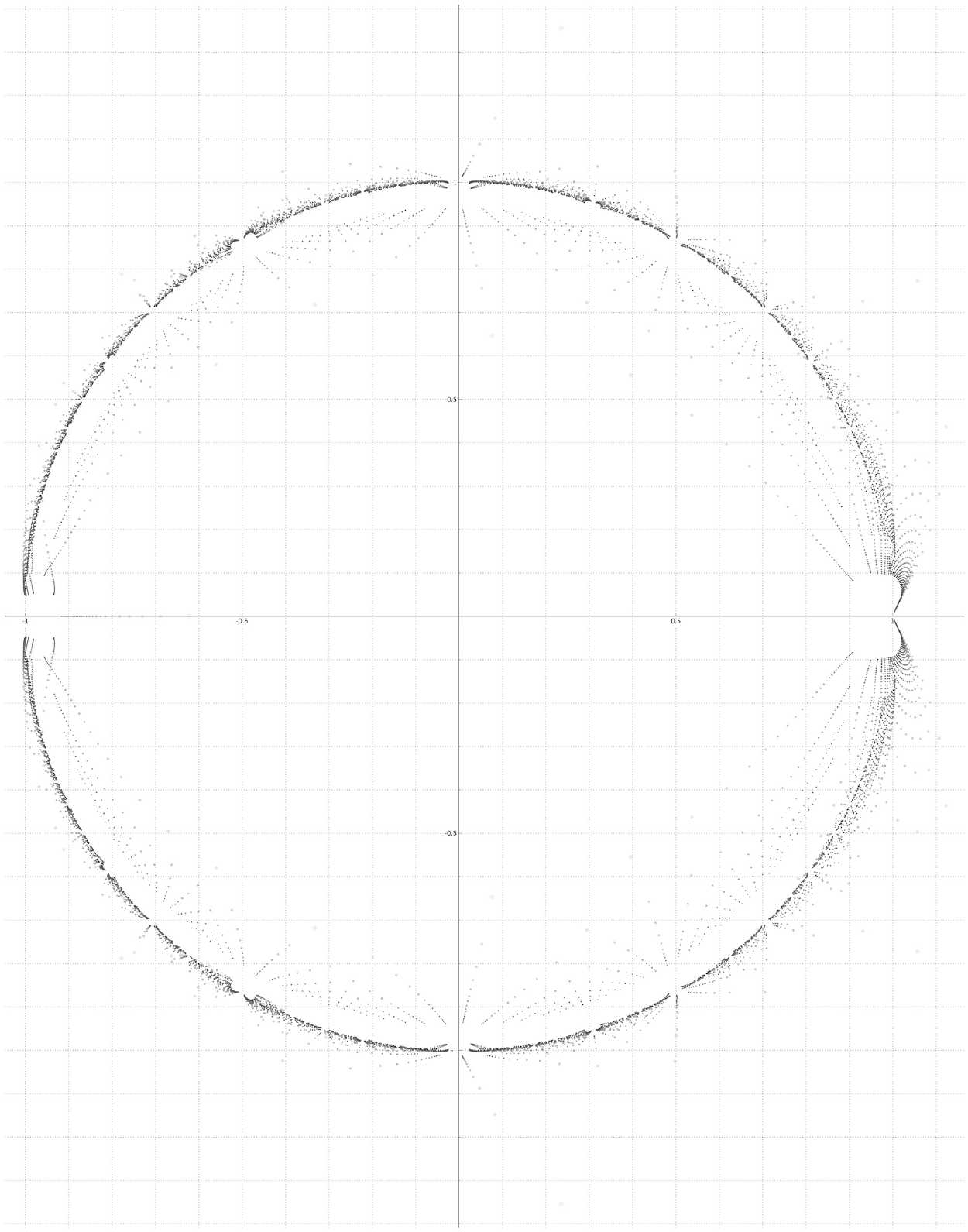}%
  \caption{ The complex zeros of $P_n(z)$ for odd $n$, where $3\le n
    \le 67$.  Darker and smaller points are used for larger $n$.
    Larger versions may be viewed at \url{http://webhome.cs.uvic.ca/~ruskey/Publications/Tatami/HoriVert.html}.
    %***WHAT DO WE DO WITH THIS URL?
    %\href{http://webhome.cs.uvic.ca/~ruskey/Publications/Tatami/HoriVert.html}{Ruskey's publications page}.
  }

 % \begin{tablenotes}
 % \item \url{http://webhome.cs.uvic.ca/~ruskey/Publications/Tatami/HoriVert.html}
 % \end{tablenotes}
\label{fig:complexPlotP}
%\end{threeparttable}
\end{figure}

The right hand side of \eref{genpn1} is the sum of two generating
functions, with odd and even powered terms, respectively.  The
sequence of coefficients of the odd power terms is $-\sum_{i=0}^{k}
2^{k-i}\binom{2i}{i}$, for $k\ge 0$ (see A082590 in \cite{Sloane}),
and that of the even power terms is $\binom{2k}{2}$, for $k\ge 1$ (see
A000984 in \cite{Sloane}).  The first few numbers $P_n(-1)$, starting
with $n=2$, are: $1$, $-1$, $ 2$, $ -4$, $ 6$, $ -14$, $ 20$, $ -48$,
$ 70$, $ -166$, $ 252$, $ -584$, $ 924$, $ -2092$, $ 3432$, $ -7616$,
$ 12870$, $ -28102$, $ 48620$, $ -104824$, $ 184756$, $ -394404$, $
705432$, $ -1494240$, $ 2704156$, $ -5692636$, $ 10400600$.
% , $ -21785872$, $ 40116600$, $ -83688344$, $
% 155117520$, $ -322494208$, $ 601080390$, $ -1246068806$, $
% 2333606220$, $ -4825743832$, $ 9075135300$, $ -18726622964$, $
% 35345263800$, $ -72798509728$, $ 137846528820$, $ -283443548276$, $
% 538257874440$, $ -1105144970992$, $ 2104098963720$, $ -4314388905704$,
% $ 8233430727600$.

\conref{Pcon}(\ref{enum:alternating}) compares the above sequence with the sum of the
absolute values of the coefficients of $P_n(z)$.  The first few of
these are listed, also starting with $n=2$: $1$, $ 3$, $ 4$, $ 10$, $
10$, $ 22$, $ 28$, $ 64$, $ 76$, $ 180$, $ 260$, $ 606$, $ 932$, $
2124$, $ 3440$, $ 7666$, $ 12872$, $ 28178$, $ 48620$, $ 104946$, $
184756$, $ 394638$, $ 705432$, $ 1494600$, $ 2704156$, $ 5693376$, $
10400600$.
% , $ 21786950$, $ 40116600$, $ 83690438$, $
% 155117520$, $ 322497282$, $ 601080390$, $ 1246075676$, $ 2333606220$, $
% 4825753816$, $ 9075135300$, $ 18726642276$, $ 35345263800$, $
% 72798537850$, $ 137846528820$, $ 283443609218$, $ 538257874440$, $
% 1105145058924$, $ 2104098963720$, $ 4314389075370$, $ 8233430727600$.

%Not referred to and removed for JCTA
% \input{Phitable.tex}

%\input{Ptable.tex}
%\clearpage

\section{Conclusions and further research}
\label{sec:conclusion}
The polynomials $P_n(z)$ exhibit numerous patterns in the signs of
their coefficients, their plots and zeros (see
\url{http://webhome.cs.uvic.ca/~ruskey/Publications/Tatami/HoriVert.html}),
their degrees, and the values of $P_n(-1)$, but yet all we have seen
is that they fall magically out of these coverings.  What is the
geometric interpretation, if there is one, for the factorisation of
$\vh_n(z)$, and how do we calculate the coefficients of $P_n(z)$
without dividing $D_n(z)$ into $\vh_n(z)$?

In the present paper we deal with $n\times n$ coverings with $n$
monominoes, and the techniques can be readily applied to $r\times c$
coverings with $r<c$ and a maximum number of monominoes.  When $m$ is
not maximum, however, we may have to deal with other features of
tatami coverings, called bidimers and vortices, described in
\cite{EricksonRuskeySchurch2011}.  The numbers of horizontal and
vertical dominoes that any given bidimer or vortex will introduce is
easily calculated, and they provide the advantage of sometimes
isolating corners of the grid, which makes diagonal flips easier to
count.  We may not, however, have the good fortune of encountering a
formula like the one in \lref{vdominoes}, which easily yields the
generating polynomial of \tref{vdominoes}.

A generating function for fixed height balanced coverings, or perhaps
some other relation between the numbers of vertical and horizontal
tiles, would extend the results on fixed height coverings in
\cite{RuskeyWoodcock2009} and \cite{EricksonRuskeySchurch2011}.

% Given the likeness of the formula
%   \begin{align*}
%     \sum_{k\ge 1} \flr{n}{2^k}^2,
%   \end{align*}
%   to the well known
%   \begin{align*}
%     \sum_{k\ge 1} \flr{n}{2^k},
%   \end{align*}
%   we were surprised not to find it elsewhere.

% \section{Acknowledgements}
%\label{sec:ack}
\textbf{Acknowledgements:} We thank Prof. Don Knuth for his valuable
communications about tatami coverings.

% http://amath.colorado.edu/documentation/LaTeX/reference/faq/bibstyles.html
\bibliographystyle{plain}
\bibliography{bibliography}

\end{document}